\documentclass[12pt]{article}
\usepackage{amsthm}
\usepackage{euscript}
\let\scr=\EuScript

\usepackage{amsfonts}
\usepackage{amsmath}
\usepackage{amssymb}
\usepackage{latexsym,color}
\usepackage{amsbsy}
\usepackage{graphicx}

\newtheorem{thm}{Theorem}[section]
\newtheorem{con}{Conjecture}[section]
\newtheorem{lemma}[thm]{Lemma}
\newtheorem{cor}[thm]{Corollary}
\newtheorem{pro}[thm]{Proposition}
\newtheorem{example}[thm]{Example}
\newtheorem{definition}[thm]{Definition}
\newtheorem{remark}[thm]{Remark}
\newtheorem{Algorithm}[thm]{Algorithm}

\newcommand{\DPn}{\mathsf{DP(n)}}
\newcommand{\EDn}{\widehat{\mathsf{DP}}{(n)}}
\newcommand{\ESn}{\widehat{\mathsf{DS}}{(n)}}
\newcommand{\OP}{\mathsf{O(P)}}
\newcommand{\Cn}{\mathsf{C(n)}}
\newcommand{\Snr}{\mathsf{S(n,r)}}
\newcommand{\Snt}{\mathsf{S(n,2)}}
\newcommand{\Sn}{\mathsf{S(n)}}
\newcommand{\DSn}{\mathsf{DS(n)}}
\newcommand{\DSmn}{\mathsf{DS(m,n)}}
\newcommand{\DPmn}{\mathsf{DP(m,n)}}
\newcommand{\Kn}{\mathsf{K(n)}}
\newcommand{\TPn}{\mathsf{TP(n)}}
\newcommand{\av}{{\cal A}}
\newcommand{\pv}{{\cal P}}
\newcommand{\Fn}{\mathsf{F(n)}}
\newcommand{\comment}[1]{}
\newcommand{\oc}{\overline{c}}
\newcommand{\od}{\overline{d}}
\newcommand{\cI}{\mathcal I}
\newcommand{\RE}{{\mathbb R}}
\newcommand{\Rl}{{\mathbb R}_{\leq}}

\def\thetheorem{\thesection.\arabic{thm}}
\def\thedefinition{\thechapter.\arabic{definition}}
\def\theremark{\thechapter.\arabic{remark}}
\def\thecorollary{\thechapter.\arabic{cor}}
\def\thelemma{\thechapter.\arabic{lemma}}
\def\theproposition{\thechapter.\arabic{pro}}
\def\theexample{\thechapter.\arabic{example}}
\def\theAlgorithm{\thechapter.\arabic{Algorithm}}

\newcommand{\ncom}{\newcommand}
\ncom{\ns}{\normalsize}
\ncom{\la}{\lambda}
\ncom{\bm}{\boldmath}
\ncom{\noi}{\noindent}
\ncom{\bq}{\begin{equation}}
\ncom{\eq}{\end{equation}}
\ncom{\beqn}{\begin{eqnarray*}}
\ncom{\eeqn}{\end{eqnarray*}}
\ncom{\ba}{\begin{array}}
\ncom{\ea}{\end{array}}
\ncom{\beq}{\begin{eqnarray}}
\ncom{\eeq}{\end{eqnarray}}
\ncom{\nno}{\nonumber}
\ncom{\hs}{\mbox{\hspace{.25cm}}}
\ncom{\rar}{\rightarrow}
\ncom{\Rar}{\Rightarrow}
\ncom{\noin}{\noindent}
\ncom{\bc}{\begin{center}}
\ncom{\ec}{\end{center}}
\ncom{\sz}{\scriptsize}
\ncom{\fpd}{\Phi(\pi^{'})}
\ncom{\fp}{\Phi(\pi) }
\ncom{\nk}{\left< \begin{array}{c}
                       n\\k \end{array} \right>}
\ncom{\nd}{1^{'},2^{'},\cdots,n^{'}}
\ncom{\R}{I\!\!R}
\ncom{\de}{\bigtriangleup (F_{2n},\leq)}
\ncom{\del}{\bigtriangleup}
\ncom{\cov}{<\!\!\!\!\cdot }
\ncom{\bt}{\begin{thm}}
\ncom{\bcon}{\begin{con}}
\ncom{\et}{\end{thm}}
\ncom{\econ}{\end{con}}
\ncom{\bl}{\begin{lemma}}
\ncom{\el}{\end{lemma}}
\ncom{\bco}{\begin{cor}}
\ncom{\ds}{\displaystyle}
\ncom{\eco}{\end{cor}}
\ncom{\bp}{\begin{pro}}
\ncom{\ep}{\end{pro}}
\ncom{\bex}{\begin{example}}
\ncom{\eex}{\end{example}}
\ncom{\bd}{\begin{definition}}
\ncom{\ed}{\end{definition}}
\ncom{\brm}{\begin{remark}}
\ncom{\erm}{\end{remark}}
\ncom{\bal}{\begin{Algorithm}}
\ncom{\eal}{\end{Algorithm}}

\ncom{\pf}{\noi {\em Proof  }}
\ncom{\be}{\begin{enumerate}}
\ncom{\ee}{\end{enumerate}}
\ncom{\s}{\subset}

\begin{document}
     
\title{\bf{\large{\textcolor{red}{\bf{The Polytope of Degree 
Partitions}}}}} 

\author{   \normalsize{\textcolor{blue}{ \bf{Amitava Bhattacharya}}} \\
         \normalsize{\em {Department of Mathematics, Statistics, and 
Computer 
Science}}\\ 
\normalsize{\em{ University of Illinois at Chicago}}\\
\normalsize{\em{ Chicago, 
Illinois 60607-7045, USA}} \\
       \textcolor{magenta}{\bf{ \small{\texttt{amitava@math.uic.edu}}}}  
\and  
\textcolor{blue}{         \normalsize{ \bf{S. Sivasubramanian and  Murali 
K. 
Srinivasan}}}
\\
    \ns{\em {Department of Mathematics}} \\
     \ns{\em{      Indian Institute of Technology, Bombay}} \\
       \ns{ \em{  Powai, Mumbai 400076, INDIA}} \\
        \textcolor{magenta}{  \bf{\small{\texttt{krishnan@math.iitb.ac.in, 
mks@math.iitb.ac.in   
}}}}}
\date{}
\maketitle

\begin{abstract}
 
{
\footnotesize{

The degree partition of a simple graph is its degree sequence 
rearranged in weakly decreasing order. 
Let $\DPn$ (respectively, $\DSn$) denote the convex hull 
of all degree partitions (respectively, degree sequences) of  
simple graphs on the vertex set $[n]=\{1,2,\ldots ,n\}$. We think of 
$\DSn$ as the symmetrization of $\DPn$ and $\DPn$ as the asymmetric part 
of $\DSn$. The polytope $\DSn$ is a well studied object (Koren {\bf 
\cite{k}}, Beissinger and Peled {\bf \cite{bp}}, Peled and Srinivasan {\bf 
\cite {ps}}, Stanley {\bf \cite{s2}}). In this paper we study the polytope 
$\DPn$ and determine its vertices (and, as a corollary, its volume), 
edges, and facets. 

Let $S(n)$ denote the set of all 2-element subsets of $[n]$. Write
elements of $S(n)$ as $(i,j)$, with $i<j$, and partially order $S(n)$ by 
componentwise $\leq$. Let $\Cn$ denote the hypercube in 
$\binom{n}{2}$-space and let $\Sn$ denote the order polytope of $S(n)$ 
(also in $\binom{n}{2}$-space). Let $M(n)$ denote the $n\times 
\binom{n}{2}$ incidence matrix of singletons vs. doubletons of $[n]$. It 
is easily seen that the image of $\Cn$ under $M(n)$ is $\DSn$. 

We define a transformation
on real sequences by repeatedly averaging over the ascending runs. Using 
this transformation we show that: 

\noi (i) The image of $\Sn$ under $M(n)$ is $\DPn$.

\noi (ii) $\DPn$ equals the Fulkerson-Hoffman-McAndrew polytope.

\noi (iii) The extreme points of $\DPn$ are the threshold partitions and
thus, for $n\geq 3$, the volume of $\DSn$ is $n! \times \mbox{ volume of 
}\DPn$.

\noi (iv) For $n\geq 4$, $\DPn$ has $2^{n-1}$ vertices, $2^{n-2}(2n - 3)$ 
edges, and $(n^2 - 3n + 12)/2$ facets.

In a different direction, we prove a hypergraph generalization of 
the well known result that a weakly 
decreasing nonnegative integral 
sequence is a degree partition of a simple graph  if and only if it is 
majorized by a 
threshold partition. 
}
}

\end{abstract}
   
{\footnotesize{     
\noi
}}     
\newpage
\begin{center}    \section{{\large  Introduction }}  
\end{center}

The degree sequence of a simple graph is a classical and 
well-studied topic in graph theory.
As explained in Chapter 3 of the book {\em Threshold graphs and 
related topics} by Mahadev and Peled {\bf \cite{mp}}, this subject goes 
hand in hand with the topic of threshold sequences, i.e., degree sequences 
of threshold graphs. Threshold sequences satisfy many of the criteria for 
degree sequences in an extremal way. In this paper we develop a  
new example of this phenomenon. Our main reference for this paper is 
Chapter 3 of the book {\bf \cite{mp}}. Another informative recent 
reference is the paper by Merris and Roby {\bf \cite{mr}}.

We consider only simple graphs. Given a simple graph $G=([n],E)$ on the 
vertex set $[n]=\{1,2,\ldots 
,n\}$, the degree $d_j$ of a vertex $j$ is the number of edges with $j$ as 
an endpoint and $d_G = (d_1,d_2,\ldots ,d_n)$ is the {\em degree sequence} 
of $G$. The {\em degree partition} of $G$ is obtained by rearranging $d_G$ 
in weakly decreasing order.
Let $DS(n)$ denote the set of all degree sequences of simple graphs 
on the vertex set $[n]$ and let $DP(n)$ denote the set of all degree 
partitions of $n$-vertex simple graphs
(note that some of the entries of a degree partition may be zero. 
It is usual to have only nonzero terms in a partition, but in this paper 
it is convenient to have this slight generality). 

Define $\DSn$, the {\em polytope of degree sequences}, to be the 
convex hull (in ${\mathbb R}^n$) of all degree sequences in $DS(n)$ and 
define $\DPn$, the {\em polytope of degree partitions}, to be the 
convex hull of all degree partitions in $DP(n)$. The study of $\DSn$ was 
begun by Koren {\bf \cite{k}} who determined its extreme points and 
showed that the linearized and symmetrized Erd\"{o}s- Gallai inequalities 
provide a
linear inequality description of $\DSn$. Beissinger and Peled {\bf \cite{bp}} 
determined the (exponential) generating function of the number of extreme 
points.  Peled and Srinivasan {\bf \cite{ps}} determined the edges and 
facets of $\DSn$ and
gave another proof of 
Koren's linear inequality 
description (we use
this proof in the present  paper).   
Finally, Stanley {\bf \cite{s2}} obtained detailed 
information on $\DSn$ including generating functions for all face numbers, 
volume, number of lattice points, and (the closely related) number of 
degree sequences (i.e., 
$\#DS(n)$). In this paper we study the polytope $\DPn$ and determine its 
vertices (and, as a corollary, its volume), edges, and facets.

Threshold graphs were introduced by Chv\'{a}tal and Hammer {\bf 
\cite{ch}} and have many different characterizations. For our purposes the 
most convenient definition is the following: a simple graph 
$G$ is {\em threshold} if every induced subgraph of $G$ has 
a dominating or an isolated vertex. Define $TS(n)$ to be 
the set of all degree sequences of threshold graphs on the vertex set 
$[n]$
and define $TP(n)$ to be the set of all degree partitions of $n$-vertex 
threshold graphs. Elements of $TS(n)$ are called {\em 
threshold sequences} and elements of $TP(n)$ are called {\em threshold 
partitions}.
If $(d_1,\ldots ,d_n)\in TP(n)$, then either $d_1 =n-1$ 
or $d_n = 0$. Using this fact inductively we easily see that 
$\#TP(n)=2^{n-1}$.

We define two further polytopes in ${\mathbb R}^n$. 
The polytope $\Kn$ is 
defined to be the solution set of the following system of linear 
inequalities:
\beq  \label{egi} 
{\ds \sum_{i\in S}} x_i \;- {\ds \sum_{i\in T}} x_i& \leq& \# S (n - 1 - 
\#T),\;\;\;\; S,T\subseteq [n],\; S\cup T\not= 
\emptyset,\; S\cap T = \emptyset. \eeq
We call $\Kn$ the {\em Koren polytope} (see {\bf \cite{k}}). Note that 
taking 
$S=\{i\},\;T=\emptyset$ gives $x_i\leq n-1$ and taking 
$S=\emptyset,\;T=\{i\}$ gives $x_i\geq 0$, showing that $\Kn$ is 
indeed a polytope. 

The polytope $\Fn$ is defined to be the solution set of the 
following system of linear inequalities:
\beq & x_1 \geq x_2 \geq \cdots \geq x_n ,& \\
\label{fhm} & {\ds \sum_{i=1}^{k}} x_i\,\; - {\ds \sum_{i=n-l+1}^n} x_i 
\leq 
k(n-1-l),\;\;\;
1\leq k+l 
\leq n.&  \eeq
We call $\Fn$ the 
{\em Fulkerson-Hoffman-McAndrew polytope} (see {\bf \cite{fhm}}). Note 
that 
(\ref{fhm}) is 
obtained 
from (\ref{egi}) by taking $S=\{1,\ldots ,k\}$ and $T=\{n-l+1,\ldots 
,n\}$. Intuitively, $\Kn$ is obtained by symmetrizing $\Fn$ and $\Fn$ is 
the asymmetric part of $\Kn$. 
Also note that 
$\Kn$ has exponentially many defining inequalities while $\Fn$ has only 
quadratically many defining inequalities. 

We now recall the Fulkerson-Hoffman-McAndrew criterion for degree 
partitions (see {\bf \cite{fhm}} and item 5 in Theorem 3.1.7 in {\bf 
\cite{mp}}). We give both the partition 
and sequence  
versions. It follows from linearizing the well-known nonlinear 
inequalities of Erd\"{o}s 
and Gallai {\bf \cite{eg}}.
\bt
\label{eg}
Let $d=(d_1,d_2,\ldots ,d_n) \in {\mathbb N}^n$. Then \\
(i) $d\in DP(n)$ if and only if $d\in\Fn$ and $d_1 + \cdots + 
d_n$ is even. \\
(ii) $d\in DS(n)$ if and only if $d\in \Kn$ and $d_1 + \cdots + 
d_n$ is even.
\et
Motivated by Theorem \ref{eg}(ii) the following result was proved in 
{\bf \cite{k,ps}}.
\bt
\label{k}
$\DSn = \Kn$, with $TS(n)$ as the set of extreme points.
\et
The main result of this paper is the following generalization and 
partition analog 
of Theorem \ref{k}.
\bt
\label{bss}
$\DPn = \Fn$, with $TP(n)$ as the set of extreme points. 
\et
We can derive most of Theorem \ref{k} as
a corollary of Theorem \ref{bss}. Given a real 
vector $x$, let $[x]$ denote the vector obtained by rearranging the 
components of $x$ in weakly decreasing order. Then it is easily seen that 
$x\in \Kn$ if and only if $[x]\in\Fn$ and using this we see that 
Theorem \ref{bss} implies that $\DSn = \Kn$ and that every extreme point 
of $\DSn$ is a threshold sequence. To complete the proof of Theorem 
\ref{k} 
we need to show that every threshold sequence is an extreme point of 
$\DSn$. This can be seen as follows. Every threshold sequence of length $n$ 
has some entry equal to $n-1$ or $0$. Using this fact inductively we see 
that no threshold sequence can be written as a convex combination of other 
degree sequences.

The argument in the preceding paragraph is not 
reversible and there is no such 
simple proof 
of Theorem \ref{bss} from  Theorem \ref{k}. Our proof of Theorem \ref{bss} 
has two main ingredients: an averaging operation on real sequences based on 
descent sets and  
Theorem \ref{k}. More precisely, we use not just the statement of 
Theorem \ref{k} but its proof from {\bf 
\cite{ps}} (for  other proofs of Theorem \ref{k}, 
see {\bf \cite{k}} and {\bf \cite{bs}}). 

Let us consider Theorem \ref{bss} from a general perspective. 
Let $P$ 
be an integral polytope in ${\mathbb R}^n$ that is closed under 
permutations of its points, i.e.,  
$x\in P$ implies $\pi . x \in P$, for all permutations $\pi$ of $[n]$. 
For example, $\DSn$ is such a polytope. Let $E$ denote the 
set of extreme points of $P$ and let $E_d \subseteq E$ denote the set 
of extreme points that have weakly decreasing coordinates. There are two 
natural ways to define the {\em asymmetric part} of $P$. In terms of 
lattice points we define the asymmetric part of $P$ as the  polytope
$$P_d = \mbox{ convex hull of }\{ (x_1,x_2,\ldots ,x_n) \in P \cap 
{\mathbb N}^n 
\;|\;x_1\geq x_2\geq \cdots \geq x_n \}.$$
In terms of linear inequalities we define the asymmetric part of $P$ as 
the polytope
$P_l$ obtained by adding the inequalities $x_1\geq 
\cdots \geq x_n$ to the list of inequalities defining $P$. 
It is easily seen that $P_d\subseteq P_l$ and $E_d \subseteq \mbox{ set of 
extreme points 
of }P_d$. Equality need not hold in these two inclusions.
For instance,
consider the polytope $P$ in ${\mathbb R}^2$ defined by: $x_1,x_2 \geq 
0,\;x_1+x_2 \leq 3$. Then it is easily checked that $P_d$ 
is strictly contained in 
$P_l$. 
If we take $P$ to be the polytope
in ${\mathbb R}^2$ defined by $x_1,x_2 \geq 
0,\;x_1+x_2 \leq 2$, then we can check  that $P_d=P_l$ but $P_d$ 
has an extreme point $(1,1)$ that is not contained  in $E_d$. 
Unexpectedly, Theorem \ref{bss} asserts that, in the case 
$P=\DSn$, we have $P_d=P_l$ and $\mbox{ set of  extreme points of 
} P_d = E_d$. Note that Theorem \ref{bss} implies that the volume of
$\DSn$ is $n!$ times the volume of $\DPn$ (for $n\geq 3$, $\DSn$ and 
$\DPn$ are full dimensional).

We now discuss another viewpoint on Theorem \ref{bss}.
Let $P$ be a finite poset. For each $p\in P$ 
introduce a variable $x_p$. The order polytope $\OP$ of $P$, 
defined by Stanley {\bf \cite{s1}}, 
is the solution set of the 
following system of linear inequalities:
\beq  \label{op}& x_p \geq x_q,\;\;p<q,\;p,q\in P,& \\
     & 0 \leq x_p \leq 1,\;\;p\in P.&  
\eeq
The constraint matrix of the inequalities (\ref{op}) is easily seen to be 
totally unimodular and thus the vertices of $\OP$ are integral. It 
follows that the vertices of $\OP$ are characteristic vectors of 
order ideals of $P$ (a subset $I\subseteq P$ is an {\em order ideal} if 
$q\in I$ and $p\leq q$ imply $p\in I$). Since the linear inequality 
description of $\OP$ is of poynomial size in $\#P$ we can optimize 
linear functions over $\OP$ in polynomial time using linear 
programming. Picard {\bf \cite{p}} showed 
that one can optimize linear functions over $\OP$ in polynomial time 
using network flows. 

Let $S(n)$ denote the set of all $2$-subsets of $[n]=\{1,\ldots ,n\}$.
We write elements of $S(n)$ as $(i,j)$, where $i < j$. 
Partially order $S(n)$ as follows: given $X=(a_1,a_2)$ and $Y=(b_1,b_2)$ 
in $S(n)$
define $X\leq Y$ if $a_i\leq 
b_i,\;i=1,2$. Let $\Sn$ denote the order polytope of 
$S(n)$ and define $\Cn$ to be 
the hypercube in $\binom{n}{2}$-space. The defining inequalities of $\Cn$ 
are (here we write the variable corresponding to a 2-element subset 
$(i,j)$ as $x_{i,j}$) 
$$ 0 \leq x_{i,j} \leq 1,\;\;(i,j)\in S(n),$$
and the defining inequalities of $\Sn$ are
\beqn  & x_{i,j} \geq 
x_{k,l}\,,\;\;(i,j)<(k,l)\in S(n),& \\
     & 0 \leq x_{i,j}\leq 1,\;\;(i,j)\in S(n).&  
\eeqn
Now let $M(n)$ denote the $n\times {\binom{n}{2}}$ incidence matrix of 
singletons vs. doubletons in $[n]$, i.e., the rows of $M(n)$ are indexed 
by $[n]$ and the columns of $M(n)$ (indexed by $S(n)$) are 
characteristic vectors of elements of $S(n)$. We think of $M(n)$ as the 
linear transformation ${\mathbb R}^{\binom{n}{2}} \rar {\mathbb R}^n$, 
$y\mapsto M(n)y$.

It is easily seen that the image of
$\Cn$ under the transformation $M(n)$ is $\DSn$. 
Theorem \ref{k} gives the defining inequalities for this image along with 
the extreme points.  

Now let us consider the image of $\Sn$ under $M(n)$. It is well 
known (see {\bf 
\cite{ch,mp}}) that the order ideals in $S(n)$ are precisely the edge sets 
of threshold graphs on the vertex set $[n]$ whose degree sequences 
$(d_1,\ldots ,d_n)$ satisfy $d_1\geq \cdots \geq d_n$. It follows that 
$M(n)(\Sn) = \TPn$, where $\TPn$ is the convex 
hull of $TP(n)$. 
As we have already seen above, no threshold 
sequence  can be written as a convex combination of other degree 
sequences. 
Thus the set of extreme points of $\TPn$ is precisely 
$TP(n)$. At this point we have the inclusions (the second of these follows 
from Theorem \ref{eg})  
$$\TPn \subseteq \DPn  \subseteq \Fn.$$
In Section 4 we prove that $\TPn = \Fn$, thereby proving 
Theorem \ref{bss}.

This paper is organized as follows. In Section 2 we recall  
two characterizations of threshold graphs. In 
Section 3 we give a simple polynomial time dynamic programming algorithm 
for optimizing 
linear functions over $\Sn$. We do not use this algorithm in the rest 
of the paper. Its main purpose is to point out that, in contrast to 
general order polytopes,  optimizing linear 
functions over 
$\Sn$ does not require linear programming 
or network flows. Since $\TPn$ is a 
linear image 
of $\Sn$ this also gives 
a polynomial time algorithm for optimizing linear functions over $\TPn$. 
In Section 4 we first introduce an averaging operation on real sequences 
based on descent sets and then use this operation to 
give another algorithm for optimizing linear 
functions over $\TPn$. 
We then show that this algorithm also optimizes 
linear functions over $\Fn$, thus showing that $\TPn = \Fn$.  In 
Section 5 
we determine the facets of $\Fn$ and give an adjacency criterion for the 
extreme points of $\Fn$. As a consequence, we obtain the following. 
\bt
For $n\geq 4$, $\Fn$ 
has  $2^{n-1}$ 
vertices, $2^{n-2}(2n-3)$ edges, and $(n^2 - 3n + 12)/2$ facets.\et  
It would be interesting to determine 
all the face numbers of $\Fn$. In particular, in analogy with the face 
numbers of the 
hypercube, we can ask whether the number of dimension $k$ faces of $\Fn$, 
for $k=0,1, \ldots ,n-1$, is of the form $P_k(n)2^{n-1-k}$, where $P_k(n)$ 
is a polynomial in $n$. 

We now consider 
the analog of $\Sn$ for hypergraphs. Let $S(n,r)$ denote the set 
of all $r$-subsets of $[n]$ (note that 
$S(n)=S(n,2)$). Write each $r$-subset of $[n]$ in increasing order of its 
elements and partially order $S(n,r)$ by componentwise $\leq$ (just as for 
$S(n)$).
An {\em $r$-graph} ( or an {\em $r$-uniform hypergraph}) is a pair 
$([n],E)$, where 
$[n]$ is the  set of {\em vertices} and $E\subseteq S(n,r)$
is the set of {\em edges}. Given an $r$-graph 
$H=([n],E)$, the degree $d_j$ of a vertex $j$ is 
the number of edges in $E$ containing $j$ and $d_H = (d_1,\ldots ,d_n)$ is 
the {\em degree sequence} of $H$. The {\em degree partition} of $H$ is 
its degree sequence rearranged in weakly decreasing order. 

An $r$-graph $([n],E)$ is an {\em $r$-ideal} if $E$ is an order ideal of 
$S(n,r)$ and is said to be 
{\em $r$-threshold} if there are 
real vertex 
weights $(c_1,c_2,\ldots ,c_n)$ such that, for $X\in S(n,r)$, we 
have $c(X)\geq 0$ if $X\in E$ and $c(X)<0$ if $X\not\in E$. The degree 
sequence (respectively, degree partition) of an $r$-threshold hypergraph 
is 
called a {\em $r$-threshold 
sequence} (respectively, {\em $r$-threshold partition}) and the degree 
sequence of an $r$-ideal is called a {\em 
$r$-ideal partition} (it is easily seen that the degree sequence of an 
$r$-ideal is a partition). It is easy to  see that
every $r$-threshold partition is an $r$-ideal partition ({\bf 
\cite{bs}}). As noted above, in the case of graphs, i.e., for $r=2$, the 
set of 
threshold and ideal partitions coincide ({\bf \cite{ch,mp}}) but for 
$r\geq 3$ there exist $r$-ideal partitions that are not $r$-threshold (an 
example with $r=3$ and $n=9$ is given in {\bf 
\cite{bs}}). 
 
Let $\Snr$ 
denote the order polytope of $S(n,r)$ and consider the image of this 
polytope under the $n\times \binom{n}{r}$ incidence matrix of singletons 
vs. $r$-subsets of $[n]$. Using the linear programming solver \texttt{ 
QSopt}, authored by David Applegate, William Cook, Sanjeeb Dash, and 
Monika Mevenkamp and available at \texttt{
http://www2.isye.gatech.edu/$\sim$wcook/qsopt}, we have checked that, 
unlike 
the case of graphs, the set of extreme points of this image is a strict 
subset of the set of all $r$-ideal partitions (our example has $n=9, 
r=3$). It is conceivable that the set of extreme 
points 
of 
the image is 
precisely the set of $r$-threshold partitions (this is true in our example 
with $n=9, r=3$) though we have no proof of this.

The convex hull of degree sequences of all $r$-graphs was studied in {\bf 
\cite{bs}} were it was shown that much of the basic theory of $\DSn$ 
generalizes to this case. In particular, the extreme points of the degree 
sequence polytope are precisely the $r$-threshold sequences. In the graph 
case the threshold sequences are extremal among all degree 
sequences  in another sense also.
It is proved in {\bf \cite{ps,rg}}
that a weakly decreasing nonnegative integral 
sequence is a degree partition (of a simple graph) if and only if it is 
majorized by a threshold partition. In Section 6 we prove a hypergraph 
version of this result replacing $r$-threshold partitions by $r$-ideal 
partitions.
\bt \label{mhd}
A weakly decreasing  nonnegative 
integral sequence is a degree partition of an $r$-graph 
if and only if it is majorized by an $r$-ideal partition.
\et
We do not know whether this characterization can be used in any way to 
give an 
efficient (i.e., 
polynomial time) 
algorithm for the long 
standing open problem of recognizing degree partitions of $r$-graphs 
(see {\bf \cite{b}}). In any case, this characterization suggests the 
problem of finding an efficient algorithm for recognizing $r$-ideal 
partitions. An efficient algorithm for recognizing $r$-threshold 
partitions 
was given in {\bf \cite{bs}} using the ellipsoid method.

\begin{center} \section{{\large Threshold graphs}} \end{center}

In this short section we recall two characterizations of threshold 
graphs. The proofs are straightforward and can be found in
{\bf \cite{ch,mp}}. Let $i\leq j$. A graph 
$T=(\{i,\ldots ,j\},E)$ on the 
vertex set $\{i,\ldots ,j\}$ is said to be a {\em proper threshold graph} 
if $T$ is 
threshold and $d_i\geq d_{i+1} \geq \cdots \geq d_j$, where $d_\ell$ is 
the 
degree of 
vertex $\ell$.

\bt \label{tg}
Let $T=([n],E)$ be a simple graph on the vertex set $[n]$. The following 
are equivalent:\\
(i) $T$ is a proper threshold graph. \\
(ii) $E$ is an order ideal in $S(n)$. \\
(iii) There exist real numbers $b_1\geq b_2 \geq \cdots \geq b_n$ such 
that $(i,j)\in E$ if and only if $b_i + b_j \geq 0$.$\Box$
\et

Consider the set $TP(n)$ of degree partitions of $n$-vertex threshold 
graphs. Partially order $TP(n)$ by componentwise $\leq$, i.e., 
$(d_1,\ldots ,d_n) \leq (e_1,\ldots ,e_n)$ iff $d_i\leq e_i$ for all $i$.
Let ${\cal O}(S(n))$ denote the poset (actually, a lattice) of all order 
ideals of $S(n)$ under 
containment.

\bt \label{lattice}
(i) $TP(n)$ is a lattice with join and meet given by componentwise 
maximum and minimum, i.e., for 
$d=(d_1,\ldots ,d_n),\;e=(e_1,\ldots ,e_n)\in TP(n)$
$$d\vee e = ({\mbox max}(d_1,e_1),\ldots ,{\mbox max}(d_n,e_n)),\;\;
d\wedge e = ({\mbox min}(d_1,e_1),\ldots ,{\mbox min}(d_n,e_n)).$$  
(ii) The map ${\cal D}: {\cal O}(S(n)) \rar TP(n)$ given by 
$${\cal D}(E) = 
\mbox{ 
degree sequence of }([n],E)$$ 
is a lattice isomorphism. $\Box$
\et

\begin{center}\section{{\large Optimizing linear functions over 
$\Sn$}}\end{center}
In this section we give a simple dynamic programming algorithm for 
optimizing linear functions over $\Sn$. 

Given real weights $c=(c_{i,j} : (i,j) \in S(n))$ consider the linear 
program
\beq \label{lp1}
&& \mbox{maximize } \ds{\sum_{(i,j) \in S(n)}} c_{i,j} x_{i,j}  \\
&& \mbox{subject to } (x_{i,j} : (i,j) \in S(n)) \in \Sn.  \nonumber
\eeq
We noted in the introduction that the extreme points of $\Sn$ are 
characteristic vectors of order ideals in $S(n)$ and thus we can solve 
(\ref{lp1}) by solving the following combinatorial optimization problem 
(where $c(I) = \sum_{(i,j) \in I} c_{i,j}$ denotes the {\em weight} of the 
order ideal $I$)
\beq \label{lp2}
&& \mbox{maximize } c(I)  \\
&& \mbox{subject to } I \in {\cal O}(S(n)).  \nonumber
\eeq

\bl \label{lui} 
Given real weights $c=(c_{i,j} : (i,j) \in S(n))$, the set of maximum 
weight order ideals is closed under union and intersection. Thus, among 
the maximum weight order ideals, there is a unique maximal and a unique 
minimal element (under containment).
\el
\pf Let $I$ and $J$ be maximum weight order ideals. Then $c(I\cup J)\leq 
c(I)$, $c(I\cap J) \leq c(J)$ and $c(I) + c(J) = c(I\cup J) + c(I\cap 
J)$. The result follows. $\Box$

\bl \label{lae}
Let $I,J \subseteq S(n)$ be order ideals such that $\chi(I)$ and 
$\chi(J)$ (the characteristic vectors of $I$ and $J$) are adjacent 
vertices of $\Sn$. Then $I\subseteq J$ or 
$J\subseteq I$.
\el
\pf There is a cost vector $c$ such that $I$ and $J$ are the only maximum 
weight ideals w.r.t $c$. The result now follows from Lemma \ref{lui}. 
$\Box$

We now give a dynamic programming algorithm to find maximum weight ideals 
in $S(n)$. Let $c=(c_{i,j} : (i,j)\in S(n))$ be a cost vector. By Theorem 
\ref{tg}(ii) order ideals in $S(n)$ are precisely edge sets of proper 
threshold graphs on $[n]$ and thus finding a maximum weight order ideal is 
equivalent to finding a proper threshold graph $T=([n],E)$ with $c(E)$ 
maximum. 
In the algorithm below, for $i \leq j$, $(\{i,\ldots ,j\}, 
E_{i,j})$ 
will be the unique 
edge maximal proper threshold graph on the vertices $\{i,\ldots ,j\}$ 
with 
maximum weight.

\vskip 1ex 
\noi {\bf Algorithm 1}\\
{\em Input}: $c=(c_{i,j} : (i,j)\in S(n)).$\\
{\em Output}: The unique edge maximal proper threshold graph on $[n]$ with 
maximum 
weight. \\
{\em Method}:\\
1. {\bf for} $i$ from $1$ to $n$ {\bf do} $E_{i,i} \leftarrow 
   \emptyset$\\
2. {\bf for} $i$ from $n-1$ downto $1$ {\bf do} \\
3. $\;\;\;\;\;${\bf for} $j$ from $i+1$ to $n$ {\bf do} \\
4. $\;\;\;\;\;\;\;\;\;\;${\bf if} ($c_{i,i+1}+c_{i,i+2} + \cdots + c_{i,j} 
                   + c(E_{i+1,j}))
                    \geq c(E_{i,j-1})$\\
5. $\;\;\;\;\;\;\;\;\;\;\;\;\;\;\;${\bf then} $E_{i,j} \leftarrow 
\{(i,i+1),(i,i+2),\ldots ,(i,j)\} \cup E_{i+1,j}$\\
6. $\;\;\;\;\;\;\;\;\;\;\;\;\;\;\;${\bf else} $E_{i,j} \leftarrow 
E_{i,j-1}$ \\
7. {\bf Output} $([n],E_{1,n})$

\bl Algorithm 1 is correct, i.e., for all $i < j$,
$(\{i,\ldots ,j\}, E_{i,j})$ is the unique 
edge maximal proper threshold graph on the vertices $\{i,\ldots ,j\}$ 
with maximum weight w.r.t. $c$.
\el
\pf By induction on $j-i$. The statement is clear for $j=i$. For $i < j$ 
consider the unique maximum weight 
edge maximal proper threshold graph $T$ on the vertices $\{i,\ldots ,j\}$ 
with edge set, say, $E$. If $i$ is dominating in $T$ then, by induction, 
we have $E = \{(i,i+1),\ldots ,(i,j)\} \cup E_{i+1,j}$. If $j$ is isolated 
in $T$ 
then, by induction, we get $E = E_{i,j-1}$. It is easy to see that $i$ is 
dominating in $T$ if and only if 
$(c_{i,i+1}+c_{i,i+2} + \cdots + c_{i,j} + 
                    c(E_{i+1,j}))
                    \geq c(E_{i,j-1})$.
That completes the proof. $\Box$

Given real numbers $c_i, i\in [n]$ consider the following linear program
\beq \label{lp3}
&& \mbox{maximize } \ds{\sum_{i \in [n]}} c_i  x_i  \\
&& \mbox{subject to } (x_i : i \in [n]) \in \TPn.  \nonumber
\eeq
We noted in the introduction that the extreme points of $\TPn$ are 
the threshold partitions in $TP(n)$ and thus we can solve 
(\ref{lp3}) by solving the following combinatorial optimization problem 
\beq \label{lp4}
&& \mbox{maximize } \sum_{i\in [n]} c_i d_i  \\
&& \mbox{subject to } (d_1,\ldots ,d_n) \in TP(n).  \nonumber
\eeq
Consider problem (\ref{lp4}). Define weights $c=(c_{i,j} : (i,j) \in 
S(n))$ 
by $c_{i,j} = c_i + c_j$. Recall the poset isomorphism ${\cal 
D}:{\cal O}(S(n))\rar 
TP(n)$ and observe that, for $d=(d_1,d_2,\ldots ,d_n)\in TP(n)$, we have 
$$\sum_{i\in [n]}c_i d_i \;=\; \sum_{(i,j)\in {\cal D}^{-1}(d)} c_{i,j} 
\;=\; 
c({\cal 
D}^{-1}(d)).$$
Thus solving (\ref{lp4}) is a special case of solving (\ref{lp2}). From 
Lemmas \ref{lui}, \ref{lae} and the poset isomorphism ${\cal D}$ we now 
have
\bl \label{ljm} 
(i) Given real weights $(c_i : i \in [n])$, the set of optimal threshold 
sequences in (\ref{lp4}) 
is closed under $\vee$ and $\wedge$. Thus, among 
the optimal threshold sequences, there is a unique maximal and a unique 
minimal element. \\
(ii) Let $d,e\in TP(n)$ be adjacent vertices of $\TPn$. Then $d$ and $e$ 
are comparable in the partial order on $TP(n)$. $\Box$
\el

In Section 5 we shall characterize comparable pairs $d,e\in TP(n)$ 
that are adjacent vertices of $\TPn$.

\begin{center} \section{\large{Repeated averaging over ascending 
runs}}\end{center}

In this section we prove Theorem \ref{bss}. We begin by 
defining an averaging operation on real sequences.

Let $c=(c_1,c_2,\ldots ,c_n) \in {\mathbb R}^n$. We define its {\em 
descent set}, denoted $Des(c)$, by
$$  Des(c) = \{ i\in [n-1] \;|\; c_i > c_{i+1} \}.$$
For instance, if 
$c=(\textcolor{red}{\bf{1,3}},\textcolor{blue}{\bf{2,7}},
\textcolor{magenta}{\bf{2,3}},\textcolor{green}{\bf{1,1,5}})$ 
then 
$Des(c)=\{2,4,6\}$. Write 
the 
descent set of $c$ as $\{i_1,i_2,\ldots ,i_k\}$, where $i_1<i_2<\cdots 
<i_k$. The subsequences 
$$c_1,c_2,\ldots ,c_{i_1}\;;\; c_{i_1+1}\cdots 
c_{i_2}\;;\;\ldots \;;\;c_{i_k+1}\cdots c_n$$
are called the {\em ascending runs} of $c$. In the example 
above 
the ascending runs are 
$\textcolor{red}{\bf{1,3}};\textcolor{blue}{\bf{2,7}};
\textcolor{magenta}{\bf{2,3}};\textcolor{green}{\bf{1,1,5}}$.

Given a real vector $c=(c_1,\ldots ,c_n)$ define $\av(c)\in {\mathbb R}^n$ 
as follows: replace each $c_i$ by the average of the elements of the 
(unique) ascending run of $c$ in which $c_i$ appears. For the example from 
the preceding paragraph we have 
$$\av(c) = 
\left(\textcolor{red}{\bf{2,2,\frac{9}{2},\frac{9}{2}}},
\textcolor{blue}{\bf{\frac{5}{2},\frac{5}{2}}},
\textcolor{magenta}{\bf{\frac{7}{3},\frac{7}{3},
\frac{7}{3}}}\right),$$ 
and
$$\av(\av(c)) = 
\left(\textcolor{red}{\bf{\frac{13}{4},\frac{13}{4},\frac{13}{4},\frac{13}{4}}},
\textcolor{blue}{\bf{\frac{5}{2},\frac{5}{2}}},
\textcolor{magenta}{\bf{\frac{7}{3},\frac{7}{3},
\frac{7}{3}}}\right).$$ 

Set $\RE^n_\geq = \{(x_1,x_2,\ldots ,x_n)\in \RE^n \;|\; x_1\geq x_2\geq 
\cdots 
\geq x_n \}$.
\bl
(i) $\av(c) = c$ if and only if  $c\in \RE^n_\geq$.\\ 
(ii) Given $c\in \RE^n$, there exists $0\leq t \leq n-1$ such that 
$\av^t(c) \in \RE^n_\geq$.
\el
\noi \pf
(i) This is clear.\\
(ii) Clearly, $Des(\av(c))\subseteq Des(c)$. If $Des(\av(c))=Des(c)$ then 
$\av(c)\in\RE^n_\geq$.  Thus each 
application of the operation $\av$ either strictly decreases the descent 
set or else the process terminates. The result follows. $\Box$

Define $\pv : \RE^n \rar \RE^n_\geq$ by $\pv(c)=\av^{n-1}(c)$. Alladi 
Subramanyam has pointed out to us that the function $\pv$ arises in the 
simply 
ordered case of isotonic regression studied in order restricted 
statistical inference (see Chapter 1 of {\bf \cite{rwd}}), where the 
following geometric interpretation is given: for $c\in \RE^n$, $\pv(c)$ is 
the unique closest point (under Euclidean distance) to $c$ in the closed, 
convex set $\RE^n_\geq$.

The next two lemmas use the function $\pv$ to reduce the problem of 
maximizing 
linear functions over $\TPn$ to that of maximizing linear functions over 
$\DSn$, where a simple greedy method works.
\bl \label{ml1}
Consider the combinatorial optimization problem (\ref{lp4})
with cost vector $c=(c_1,\ldots ,c_n)\in \RE^n$.
 
\noi (i) Suppose that $c_i \leq c_{i+1}$ for some 
$i \leq n-1$. Let $d^* = (d_1^*,d_2^*,\ldots ,d_n^*)$ be the unique 
maximal optimal solution to (\ref{lp4}). 
Then $d_i^* = d_{i+1}^*$.

\noi (ii) Suppose that $c_i \leq c_{i+1}$ for some 
$i \leq n-1$. Let $d^* = (d_1^*,d_2^*,\ldots ,d_n^*)$ be the unique 
minimal optimal solution to (\ref{lp4}). 
Then $d_i^* = d_{i+1}^*$.

\noi (iii) Suppose that $c_i < c_{i+1}$ for some 
$i \leq n-1$. Let $d^* = (d_1^*,d_2^*,\ldots ,d_n^*)$ be any 
optimal solution to (\ref{lp4}). 
Then $d_i^* = d_{i+1}^*$.
\el
\noi \pf We prove part (i). The proofs for parts (ii) and (iii) are 
similar.

The proof is by induction on $n$, the case $n=2$ being clear. Let $n\geq 
3$ 
and consider the following three cases:

\noi (a) $2\leq i <  i+1 \leq n-1$: Either $d_1^* = n-1$ or $d_n^* = 0$. 
In 
the 
first case $(d_2^* - 1, \ldots ,d_n^* -1)$ is the unique maximal optimal 
solution to 
(\ref{lp4}) with cost vector $(c_2^*,\ldots ,c_n^*)$ and in the second 
case $(d_1^*,\ldots , d_{n-1}^*)$ is the unique maximal optimal solution 
to (\ref{lp4}) 
with cost vector $(c_1^*,\ldots ,c_{n-1}^*)$. By induction we now see that
$d_i^* = d_{i+1}^*$.

\noi (b) $i=1, i+1 =2$: Let $T=([n],E)$ be the proper threshold graph with 
degree 
sequence $d^*$, i.e., $E= {\cal D}^{-1}(d^*)$ and assume that $d_1^* > 
d_2^*$. Since $E$ is an order 
ideal of $S(n)$ we see that,
for some $2\leq j  < l$, the vertices adjacent to $1$ are 
$\{2,3,\ldots ,l\}$ and the vertices adjacent to 
$2$ are $\{1,2, \ldots ,j\} - \{2\}$. Let $E' = \{(1,k) \;|\; 
j < k \leq l\}$ and $E'' = \{(2,k) \;|\; j < k \leq l\}$. Note that 
$T' = ([n],E-E')$ is a proper threshold graph and thus, since $d^*$ 
is an optimal solution to (\ref{lp4}), it follows that $c(E')\geq 0$ (for 
a subset $X\subseteq S(n)$ we set $c(X) = \sum_{(i,j)\in X} c_i + c_j$). 
Since $c_2\geq c_1$ we have $c(E'')\geq 0$ and thus, since $T'' 
=([n],E\cup E'')$ is a proper threshold graph, the degree 
sequence of $T''$ is also an optimal solution to (\ref{lp4}), 
contradicting the maximality of $d^*$. So $d_1^* = d_2^*$.

\noi (c) $i=n-1, i+1=n$: Similar to case (b).$\Box$

\bl \label{ml2}
Let $c\in\RE^n$. Consider two instances of the combinatorial optimization 
problem (\ref{lp4}), one with cost vector $c$ and another with cost vector 
$\pv(c)$. Then 

\noi (i) The unique maximal optimal 
solutions to these two instances are equal.

\noi (ii) The unique minimal optimal 
solutions to these two instances are equal.

\el
\pf We prove part (i). The proof for part (ii) is similar.
We show that the unique maximal optimal solutions to (\ref{lp4}) with 
cost vectors $c$ and $\av(c)$ are the same. This will prove part (i). 

Let $d^*=(d_1^*,\ldots ,d_n^*)$ be the unique maximal optimal solution to 
(\ref{lp4}) 
with cost vector $c$ and let $e^*=(e_1^*,\ldots ,e_n^*)$ be the 
unique maximal optimal 
solution to (\ref{lp4}) with cost vector $\av(c)$.

Let $c=(c_1,\ldots ,c_n)$ and $\av(c)=(b_1,\ldots ,b_n)$. Write
$Des(c)=\{i_1,i_2,\ldots ,i_k\}$, where $i_1 < i_2 < \cdots < i_k$. Put 
$i_0 
= 0, i_{k+1}=n$ and, for $\ell=1,2,\ldots ,k+1$, set
\beq \label{bt} B_\ell = \{ i_{\ell-1} + 1,i_{\ell-1} + 2 , \ldots ,i_\ell 
\},\eeq 
i.e., $B_\ell$ is the set of indices of the $\ell\,$th ascending run of 
$c$.

By Lemma \ref{ml1} and the definition of the map $\av$ we have
$$d_i^* = d_j^* \mbox{ and } e_i^* = e_j^* \mbox{ whenever } i,j\in 
B_\ell, 
\mbox{ for some } \ell.$$

We now have, using the definition of the map $\av$,
\beqn
c_1d_1^* + c_2d_2^* + \cdots + c_nd_n^* & = & 
\sum_{\ell=1}^{k+1} \left\{ 
\sum_{s\in B_\ell} c_s\right\} d_{i_\ell}^* \\
   & = & \sum_{\ell=1}^{k+1} \left\{ 
\sum_{s\in B_\ell} b_s\right\} d_{i_\ell}^* \\
 & = & b_1d_1^* + b_2d_2^* + \cdots + b_nd_n^*.
\eeqn
Similarly we can show
\beqn
c_1e_1^* + c_2e_2^* + \cdots + c_ne_n^* & = & 
b_1e_1^* + b_2e_2^* + \cdots + b_ne_n^*.
\eeqn
    
Now we use the fact that $d^*$ is optimal for the cost vector $c$ and 
$e^*$ is optimal for the cost vector $\av(c)$. We have
$$\sum_{i=1}^n c_id_i^* \geq
\sum_{i=1}^n c_ie_i^* =
\sum_{i=1}^n b_ie_i^* \geq
\sum_{i=1}^n b_id_i^* =
\sum_{i=1}^n c_id_i^*.$$
It follows that 
$\sum_{i=1}^n c_id_i^* =
\sum_{i=1}^n c_ie_i^*$ and  
$\sum_{i=1}^n b_id_i^* =
\sum_{i=1}^n b_ie_i^*.$ 

Since $d^*$ is the unique maximal solution to (\ref{lp4}) with cost vector 
$c$ we have $e^* \leq d^*$ 
and since $e^*$ is the unique maximal solution to (\ref{lp4}) with cost 
vector 
$\av(c)$ we have $d^* \leq e^*$. Thus $d^* = e^*$. $\Box$

We can now give our second algorithm to solve the optimization problem 
(\ref{lp4}).

\vskip 1ex 
\noi {\bf Algorithm 2}\\
{\em Input}: $c=(c_1,\ldots,c_n)\in {\mathbb R}^n.$\\
{\em Output}: The unique maximal $d^*=(d_1^*,\ldots ,d_n^*)\in TP(n)$ 
maximizing 
$\sum_{i=1}^n c_id_i$ over all $(d_1,\ldots ,d_n)\in TP(n).$\\
{\em Method}:\\
1. $(b_1,\ldots ,b_n)\leftarrow \pv(c)$\\
2. $E\leftarrow \emptyset$\\
3. {\bf for all} ($(i,j)\in S(n)$ with $b_i + b_j \geq 0$) {\bf do}
           $E \leftarrow E \cup \{(i,j)\}$ \\
4. $(d_1^*,\ldots ,d_n^*)\leftarrow \mbox{ degree sequence of } ([n],E)$

\bl \label{ac} Algorithm 2 is correct.
\el
\pf Consider the problem of maximizing $b_1d_1 +\cdots + b_nd_n$ over all 
degree sequences $(d_1,\ldots ,d_n)\in DS(n)$. An edge $(i,j)$ 
contributes $b_i + b_j$ to the objective function and hence it follows 
that among all simple graphs on $[n]$ whose degree 
sequences maximize $\sum_{i=1}^nb_id_i$, the graph $([n],E)$ 
computed in Step 3 of Algorithm 2 is the unique maximal one (under 
containment of edges).
Since $b_1\geq b_2 \geq \cdots \geq b_n$, it follows from 
Theorem \ref{tg}(iii) that $d^*$ is a threshold 
partition. Using Lemma 
\ref{ml2} we  now see that $d^*$ is the unique maximal optimal solution to 
(\ref{lp4}) with cost vector $c$. $\Box$

\noi {\bf Remarks} (i) To get the unique minimal optimal solution in 
Algorithm 
2 we 
change line 3 to the following:\\
3. {\bf for all} ($(i,j)\in S(n)$ with $b_i + b_j > 0$) {\bf do}
           $E \leftarrow E \cup \{(i,j)\}$ 

\noi (ii) As observed in the introduction it is easy to see that every 
element 
of $TS(n)$ (respectively, $TP(n)$) is an 
extreme point of $\DSn$ (respectively, $\DPn$). Algorithm 2 shows that the 
converse statement for $\DSn$, namely, that 
every extreme point of $\DSn$ is in $TS(n)$ also has a short direct proof 
not involving the polytope $\Kn$. 
In contrast, our proof that every extreme point of $\DPn$ is in $TP(n)$ is 
indirect and uses the polytope $\Fn$.

\vskip 1ex
We can now give the proof of our main theorem.

\noi
\pf  (of Theorem \ref{bss}) Let $c=(c_1,\ldots ,c_n) \in \RE^n$ and let 
$d^*=(d_1^*,\ldots ,d_n^*)$ be the output of Algorithm 2. We show that 
$d^*$ maximizes $\sum_{i=1}^n c_ix_i$ over $\Fn$. This will prove that 
$\TPn = \Fn$ and since, as already observed in the introduction, the set 
of extreme points of $\TPn$ is $TP(n)$, this also proves Theorem 
\ref{bss}.

Write the constraints of $\Fn$ in the form 
\beq \label {rv} & -x_i + x_{i+1} \leq 0,\;\;\;1\leq i\leq n-1.& \\
\label{fhmc} & {\ds \sum_{i=1}^{k}} x_i\,\; - {\ds \sum_{i=n-l+1}^n} x_i 
\leq 
k(n-1-l),\;\;\;
1\leq k+l 
\leq n.&  \eeq
We shall show that the row vector $c$ can be written as a nonnegative 
rational combination of the row vectors of lhs coefficients of 
those constraints from (\ref{rv}) and (\ref{fhmc}) that are satisfied with 
equality by $d^*$. By the (weak) duality theorem of linear 
programming, this will prove the result.

Write $Des(c)=\{i_1,i_2,\ldots ,i_k\}$, where $i_1 < i_2 < \cdots < i_k$. 
For $\ell=1,\ldots ,k+1$, define $B_\ell$ to be the set of indices of
the $\ell\,$th 
ascending run of $c$, as in (\ref{bt}). Let $\pv(c) = 
(b_1,\ldots ,b_n)$. Since the ascending runs of ${\cal A}^j(c)$ are 
contained in the ascending runs of ${\cal A}^{j+1}(c)$ it follows from the 
definition of $\pv(c)$ that $b_i = b_{i+1}$ whenever $i,i+1 \in B_\ell$, 
for 
some $\ell$. It now follows from Algorithm 2 that the inequality $-x_i + 
x_{i+1} \leq 0$ is satisfied with equality by $d^*$ whenever $i,i+1\in 
B_\ell$, for some $\ell$.

For $1\leq i\leq n-1$ define
$v_i = (0,\ldots ,0,-1,1,0,\ldots ,0),$
with $-1$ in the $i\,$th spot and $1$ in the $(i+1)\,$st spot, i.e., 
$v_i$ is 
the row vector of lhs coefficients of the $i\,$th inequality from 
(\ref{rv}).

For $\ell=1,\ldots ,k+1$, let $m_\ell = \mbox{min }B_\ell = i_{\ell-1} +
1, M_\ell = \mbox{max }B_\ell = i_\ell$ 
and let 
$a_\ell = \frac{\sum_{i\in B_\ell} c_i}{  \#B_\ell}$ 
be the 
average of the elements of the $\ell\,$th ascending run of $c$. 
We assert 
that
\beq\label{ci}
c = \av(c) + \sum_{\ell=1}^{k+1} \sum_{i=m_\ell}^{M_\ell - 
1} \left\{ (i+1- m_\ell)a_\ell - (c_{m_\ell}+\cdots 
+c_i)\right\} v_i.
\eeq

Before proving the assertion we note that the coefficient of $v_i$ in 
(\ref{ci}) is nonnegative since $a_\ell$ is the average of  
the $\ell\,$th ascending run and the entries in this run are 
weakly 
increasing. Also note 
that for every $v_i$ that appears in (\ref{ci}) the corresponding 
inequality in (\ref{rv}) is satisfied with equality by $d^*$. 

Let $1\leq \ell \leq k+1$. We now show that, for $i\in B_\ell$, the 
$i\,$th 
coordinate on the rhs of (\ref{ci}) is $c_i$. This will prove 
(\ref{ci}). The $m_\ell\,$th coordinate on the rhs of (\ref{ci}) is 
$a_\ell 
- (a_\ell 
- c_{m_\ell}) = c_{m_\ell}$ and the $M_\ell\,$th coordinate on the rhs of 
(\ref{ci}) 
is $a_\ell + ((M_\ell - m_\ell)a_\ell - (c_{m_\ell} + \cdots + c_{M_\ell - 
1})) = c_{M_\ell}$.
For $m_\ell < i < M_\ell$, the $i\,$th coordinate on the rhs of (\ref{ci}) 
is
$a_\ell + ((i - m_\ell)a_\ell - (c_{m_\ell} + \cdots + c_{i-1})) - 
((i+1-m_\ell)a_\ell - 
(c_{m_\ell} + \cdots + c_i)) = c_i$.

Repeating the above procedure for $\av(c), {\cal A}^2(c), \ldots $ we 
eventually arrive at an expression of the form 
$$ c = \pv(c) + \sum \alpha_i v_i,$$
where the sum is over all $1\leq i \leq n-1 $ with $d_i^* = d_{i+1}^*$ and 
$\alpha_i$ is a 
nonnegative rational for all $i$.

We observed in the proof of Lemma \ref{ac} that $d^*$ is an optimal degree 
sequence for the cost vector $\pv(c)$. Therefore, by Theorem \ref{k}, 
$\pv(c)$ can be written as a 
nonnegative rational combination of the row vectors of lhs coefficients 
of those inequalities from (\ref{egi}) 
that are satisfied with equality by $d^*$. In fact, the proof of Theorem 
4.2 from {\bf \cite{ps}} (also see the proof of Theorem 3.3.14 from {\bf 
\cite{mp}}) 
shows how 
to write $\pv(c)$ as a nonnegative 
rational combination of the row vectors of lhs coefficients of those 
inequalities from (\ref{fhmc}) that are 
satisfied with equality by $d^*$. 

That completes 
the proof of Theorem \ref{bss}. $\Box$

\begin{center}\section{{\large  Facets and edges of  $\DPn$ }}  
\end{center}

In this section we determine the facets and  edges of the 
polytope $\DPn$. Most of the steps needed to determine the facets of 
$\DPn$ 
are similar to the case of $\DSn$ and therefore we just quote the 
corresponding results from {\bf \cite{mp}}. 
%We make an exception in the 
%case of Lemma \ref{eq} below, where we reproduce the proof from {\bf 
%\cite{mp}} as this result is important for determining the edges of 
%$\DPn$ 
%(which is different from the $\DSn$ case). 

In order to identify the facets of $\DPn$ we need to know 
its dimension and the dimension of a related polytope. For $m,n\geq 1$, 
let $K(m,n)$ denote the complete bipartite graph with bipartition 
$\{1,\ldots ,m\}$ and $\{m+1,\ldots ,n\}$. For a spanning subgraph of 
$K(m,n)$ we define its {\em degree bipartition} as the sequence 
obtained by rearranging the first $m$ and last $n$ components of its 
degree sequence 
$(d_1,\ldots ,d_m,d_{m+1},\ldots ,d_n)$ into 
weakly decreasing order. 
Define $\DSmn$ (respectively, $\DPmn$) to be the convex 
hull of degree sequences (respectively, degree bipartitions) of spanning 
subgraphs of $K(m,n)$.

\bl \label{dil}
(i) $\mbox{dim }{\mathsf{DP(1)}} = 0$, 
        $\mbox{dim }{\mathsf{DP(2)}} = 1$.\\
(ii) $\mbox{dim }{\mathsf{DP(m,n)}} = m+n-1$, for $m,n\geq 1$.\\
(iii) $\mbox{dim }{\mathsf{DP(n)}} = n$, for $n\geq 3$. 
\el
\pf (i) Clear.\\
(ii) and (iii) This follows from Lemma 3.3.16 in {\bf \cite{mp}} which 
proves  
$\mbox{dim }{\mathsf{DS(m,n)}} = m+n-1$, for $m,n\geq 1$ and
$\mbox{dim }{\mathsf{DS(n)}} = n$, for $n\geq 3$. 
A direct proof can also be given.
$\Box$
%
%\bl \label{eq} 
%Let $G=([n],E)$ be a simple graph on the vertex set $[n]$ whose degree 
%sequence $d=(d_1, d_2, \ldots ,d_n)$ satisfies 
%\beq
%& {\ds \sum_{i=1}^{k}} d_i\,\; - {\ds \sum_{i=n-l+1}^n} d_i =
%k(n-1-l),\eeq
%for some $1\leq k+l \leq n$. Then \\
%(i) $\{1,\ldots ,k\}$ induces a clique in $G$.\\
%(ii) $\{n-l+1,\ldots ,n\}$ induces a stable set  in $G$.\\
%(iii) Every vertex in $\{k+1, \ldots , n-l\}$ is adjacent to every vertex 
%in $\{1, \ldots ,k\}$.\\
%(iv) No vertex in $\{k+1, \ldots , n-l\}$ is adjacent to any  vertex 
%in $\{n-l+1, \ldots ,n\}$.
%\el
%\pf This is Lemma 3.3.13 from {\bf \cite{mp}}. $\Box$ 
%Let $C\subseteq E$ denote the edges of $G$ with one endpoint in 
%$S=\{1,\ldots ,k\}$ and 
%another in $T=\{n-l+1,\ldots ,n\}$. The degree of a vertex $i\in S$ in 
%the 
%induced subgraph of 
%$G$ on $[n]-T$ is at most $n-1- l$ and thus
%\beq  \label{egi1} 
%\sum_{i=1}^k d_i &\leq&  k (n - 1 - 
%l) + \#C. \eeq
%
%The set of edges of $G$ incident with a vertex of $T$ includes $C$ and 
%hence
%\beq  \label{egi2} 
%-\sum_{i=n-l+1}^n d_i &\leq& -\;\#C. \eeq
%
%Adding the two inequalities above we see that
%\beq  \label{egi3} 
%\sum_{i=1}^k d_i \;- \sum_{i=n-l_1}^n d_i& \leq& k(n - 1 - 
%l). \eeq
%
%It is easily seen that equality holds in (\ref{egi1})
%if and only if conditions (i) and (iii) of the lemma 
%hold and equality holds in (\ref{egi2})
%if and only if conditions (ii) and (iv) of the lemma 
%hold. $\Box$

\bt \label{facets}
For $n\geq 4$ the facets of $\DPn$ are given by \\
(i) $x_i \geq x_{i+1},\;\;i=1,\ldots ,n-1.$\\
(ii) ${\ds \sum_{i=1}^{k}} x_i\,\; - {\ds \sum_{i=n-l+1}^n} x_i 
\leq 
k(n-1-l),$ \\
where $k=1, l=0 \mbox{ or } k=0, l=1 \mbox{ or }k,l\not= 0,\;\;
k+l = 2,3,\ldots, n-3,n.$\\
It follows that, for $n\geq 4$, $\DPn$ has $\ds{\frac{n(n-3)}{2}} + 6$ 
facets. 
\et
\pf By Theorem \ref{bss} and the full dimensionality of $\DPn$ it follows 
that every facet of $\DPn$ is of the form $x_i \geq x_{i+1}$, for some 
$i=1,\ldots ,n-1$ or of the form (\ref{fhm}), for some $1\leq k+l \leq n$.

We first show that $x_i \geq x_{i+1},\;i=1,\ldots ,n-1$ determines a 
facet for $n\geq 3$. Since $(0,\ldots ,0)$ satisfies $x_1=\cdots = x_n$ it 
is enough to produce $n-1$ linearly independent elements of $TP(n)$ 
satisfying $x_i = x_{i+1}$, for $i=1,\ldots ,n-1$. This we do by induction 
on $n$. For $n=3$, $\{(1,1,0), (2,2,2)\}$ and $\{(2,1,1), (2,2,2)\}$ are 
sets of 2 linearly independent elements in $TP(3)$ satisfying $x_1=x_2$ 
and $x_2=x_3$ respectively.

Let $n\geq 4$. We consider two cases.

(a) $1\leq i \leq n-2$: 
By induction, there is a set $S$ 
of $n-2$ linearly independent elements of $TP(n-1)$ satisfying $x_i = 
x_{i+1}$. 
Let $S' \subseteq TP(n)$ 
denote the set obtained from $S$ by adding a zero at the end of every 
element of $S$. Then it is easily checked that $S' \cup \{(n-1,\ldots 
,n-1)\}$ is a set of $n-1$ linearly independent elements of $TP(n)$ 
satisfying $x_i = x_{i+1}$. 

(b) $i=n-1$: By induction, there is a set $S$ 
of $n-2$ linearly independent elements of $TP(n-1)$ satisfying $x_{n-2} = 
x_{n-1}$. 
Let $S''$ 
denote the set obtained from $S$ by adding a zero at the beginning of 
every 
element of $S$ and let $\Delta = (n-1,1,\ldots,1)\in TP(n)$. 
Then it is easily checked that $\{\Delta + u \;:\; u\in S''\} \cup 
\{\Delta\}$
is a set of $n-1$ linearly independent elements of $TP(n)$ 
satisfying $x_{n-1} = x_n$. 

Now we have to show that if $n\geq 4$ and  (\ref{fhm}) is a facet of 
$\DPn$ then 
$k,l$ satisfy the conditions listed under (ii) of the statement. This can 
be done using Lemma \ref{dil} in the same way as 
Theorem 3.3.17 (determining the facets of $\DSn$) is deduced from Lemma 
3.3.16 in {\bf 
\cite{mp}} and we  omit the 
proof. 

For $n\geq 4$, the number of pairs $(k,l)$ of positive integers with 
$k+l = p$ ($2\leq p \leq n$) is $p-1$ and thus
the number of facets of $\DPn$ is given by
$$n - 1 + 2 + \left(\sum_{p=2}^n (p-1)\right) - (n-3) - (n-2) = 
\frac{n(n-3)}{2} + 
6.\;\;\; \Box$$

By an  {\em interval} of $[n]$ we mean a nonempty 
subset $I\subseteq [n]$ of the form $I=\{i,i+1,\ldots ,j\}$, for some $i 
\leq j$. An interval $I$ is {\em nontrivial} if $\#I \geq 2$. 
Given a nontrivial interval $I\subseteq [n]$, define the vector 
$\Omega_I = (a_1,\ldots ,a_n)$ by $a_i = 0$ if $i\not\in I$, and $a_i = 
\#I - 1$ if $i\in I$.

Consider a disjoint pair of intervals of $[n]$. Write this disjoint pair 
as $(I,J)$, where $\mbox{max }I < \mbox{ min }J$. Define the vector 
$\Omega_{I,J} = (a_1,\ldots ,a_n)$ by $a_i = 0$ if $i\not\in I \cup J$,
$a_i = \#J$ if $i\in I$, and $a_i = \#I$, if $i\in J$.

\bt \label{edges}
Let $n\geq 3$ and let $d,e \in TP(n)$. Then $d$ and $e$ are adjacent 
extreme points of 
$\DPn$ if and only if $d$ and $e$ are comparable {\em (}in the partial 
order on 
threshold partitions{\em )} and 
$$|d-e| = \Omega_I \mbox{ or } \Omega_{I,J}\;,$$ 
for some 
nontrivial interval $I$ or some disjoint pair of intervals $(I,J)$.
\et
\pf Let $d=(d_1,\ldots ,d_n)$ and $e=(e_1,\ldots ,e_n)$.

(if) Assume $e\leq d$. The proof is by induction on $n$. The case $n=3$ is 
easily verified (there are four threshold partitions of length 3 and any 
two of them are adjacent). Let $n\geq 4$ and consider the following three 
cases.

(a) $d_n = e_n = 0$: By 
induction, $(d_1,\ldots ,d_{n-1})$ and $(e_1,\ldots 
,e_{n-1})$ are adjacent extreme points of $\mathsf{DP(n-1)}$ ($n$ cannot 
be a member of $I$ or $J$). The facet 
of $\DPn$ given by $x_n \geq 0$ is isomorphic to $\mathsf{DP(n-1)}$ and 
since 
$d$ and 
$e$ lie on this facet it follows that they are adjacent.

(b) $d_1 = e_1 = n-1$: By induction $(d_2 - 1,\ldots ,d_{n} - 1)$ and 
$(e_2 - 1,\ldots 
,e_{n}-1)$ are adjacent extreme points of $\mathsf{DP(n-1)}$. The facet 
of $\DPn$ given by $x_1 \leq  n-1$ is isomorphic to $\mathsf{DP(n-1)}$ and 
since $d$ 
and 
$e$ lie on this facet it follows that they are adjacent.

(c) $d_1 = n-1$ and $e_n =  0$: In this case $d-e$ will have nonzero first 
and last components. We consider two subcases.

(i) $d-e = \Omega_{[n]}$: We must have $e=(0,\ldots ,0)$ and 
$d=(n-1,\ldots , n-1)$. Since $d$ and $e$ both satisfy the $n-1$ linearly 
independent defining inequalities $x_i \geq x_{i+1}$, $i=1,\ldots ,n-1$ 
with equality they are 
adjacent.

(ii) $d-e = \Omega_{I,J}$ for some disjoint intervals $(I,J)$ with $1\in 
I$ and $n\in J$: Let $T_d$ and $T_e$ be the proper threshold graphs on the 
vertex set $[n]$ with degree sequences $d$ and $e$ respectively.
Write $d-e = (r,\ldots ,r,0,\ldots ,0,s,\ldots ,s)$ where 
$r=\#J$, $s=\#I$, $r$ is repeated $s$ times and $s$ is repeated $r$ times. 
Now $e_1 = d_1 - r = n-1-r $ and hence $e_{n-r+1}=\cdots = e_{n}=0$. 
Similarly we can show that $d_1=\cdots =d_s = n-1$. 

We now have
\beqn
d&=&(n-1,\ldots ,n-1,\ldots ,s,\ldots ,s),\\
e&=&(n-1-r,\ldots ,n-1-r,\ldots ,0,\ldots ,0),\\
d_i&=&e_i,\;\;\;i=s+1,\ldots ,n-r,
\eeqn
where $n-1$ is repeated $s$ times and $s$ is repeated $r$ times (in $d$) 
and
$n-1-r$ is repeated $s$ times and $0$ is repeated $r$ times (in $e$).

It follows that 
\begin{itemize}
\item the vertices $\{1,\ldots ,s\}$ induce a clique in both 
$T_d$ and $T_e$,
\item the vertices $\{n-r+1,\ldots ,n\}$ induce a stable set  in both 
$T_d$ and $T_e$, 
\item  every vertex of $\{s+1,\ldots ,n-r\}$ is connected to 
every vertex in $\{1,\ldots ,s\}$ in both $T_d$ and $T_e$, 
\item no vertex in $\{s+1,\ldots ,n-r\}$ is 
connected to any vertex in $\{n-r+1,\ldots ,n\}$ in both $T_d$ and $T_e$.
\end{itemize}
Thus the proper threshold subgraphs of $T_d$ 
and $T_e$ induced on the vertices 
$\{s+1,\ldots ,n-r\}$ have identical degree sequences and hence, by 
Theorem \ref{lattice}, are identical, say $T$. By Theorem \ref{bss} there 
exist real coefficients $(c_{s+1},\ldots ,c_{n-r})$ such that $T$ is the 
unique proper threshold graph on $\{s+1,\ldots ,n-r\}$ of maximum weight 
(i.e., whose degree sequence maximizes $\sum_{i=s+1}^{n-r} c_i x_i$).

Let $M = \mbox{max }\{|c_i|\;:\;s+1\leq i\leq n-r\}$. Choose
positive numbers $0 < c_1 < c_2 < \cdots < c_s$ and 
negative numbers 
$c_{n-r+1} < \cdots < c_n < 0$ with $c_1 > M$, $|c_n| > M$, and 
$r(c_1 + \cdots + c_s) = -s(c_{n-r+1} + 
\cdots + c_n)$. Put $c=(c_1,\ldots ,c_n)$ and consider problem 
(\ref{lp4}). It is easy to check, using Lemma \ref{ml1}(iii) and the 
uniqueness of $T$, that 
an optimal threshold partition with 
last component 0 must be $e$ and an optimal threshold partition with first 
component $n-1$ must be $d$. Since $\sum_{i=1}^n c_i d_i = \sum_{i=1}^n 
c_i e_i$
it follows that $d$ and $e$ are the only optimal threshold 
partitions and thus are adjacent.

\vskip 1ex
\noi
(only if) By Theorem \ref{ljm}, $e$ and $d$ are comparable, say $e\leq d$. 
There is a real cost vector 
$c\in \RE^n$ such that $d$ and $e$ are the only two 
optimal solutions to problem (\ref{lp4}). Then $e$ must be the unique 
minimal solution and $d$ must be the unique maximal solution. Thus 
Algorithm 2 will produce $d$ as output and Algorithm 2, with step 3 
modified as in remark (i) following the algorithm, will produce $e$ as 
output.

Let $\pv(c)=(b_1,b_2,\ldots ,b_n)$ and 
write $Des(\pv(c))=\{i_1,i_2,\ldots ,i_k\}$, where $i_1 < i_2 < \cdots < 
i_k$. 
For $\ell=1,\ldots ,k+1$, define $B_\ell$ to be the set of indices of 
the $\ell\,$th ascending run of $\pv(c)$, as in (\ref{bt}). We 
have 

(i) $b_i = b_j$ whenever $i,j\in B_\ell$, for some $\ell$. For 
$\ell=1,\ldots 
,k+1$, define $a_\ell = b_i$, for (any) $i\in B_\ell$. Note that $a_1 > 
a_2 > 
\cdots > a_{k+1}$. 

(ii) Set $\Sigma_0 = \{ \ell\;:\; 1\leq \ell \leq k+1,\;a_\ell=0\}$. Note 
that 
$\#\Sigma_0 \leq 1$. 

(iii) Set $\Sigma_c = \{ (i,j)\;:\; 1\leq i < j \leq k+1,\;a_i + a_j = 
0\}$. Note that the ordered pairs in $\Sigma_c$ are disjoint and 
incomparable (in the partial order on $S(k+1)$).

Let $E_d$ and $E_e$ be the edge sets of the unique proper threshold graphs 
on $[n]$ with degree sequences $d$ and $e$ respectively. Then from 
Algorithm 2 (and the subsequent remark (i)) we see that
$$E_d - E_e = \left( \biguplus_{i\in \Sigma_0} \binom{B_i}{2}\right) 
      \biguplus \left( \biguplus_{(i,j)\in \Sigma_c} B_i \times B_j 
\right),$$
where $\biguplus$ denotes disjoint union and $\binom{X}{2}$ (for $X$ a set 
of integers) denotes the set 
of all ordered pairs $(i,j)$ with $i<j$ and $i,j\in X$.

Now we observe that, for $i\in \Sigma_0$, $E_e \cup \binom{B_i}{2}$ and 
$E_d - \binom{B_i}{2}$ are both order ideals in $S(n)$ and, for $(i,j)\in 
\Sigma_c$, $E_e \cup (B_i\times B_j)$ and
$E_d - (B_i\times B_j)$ are also both order ideals in $S(n)$. Since $d$ 
and 
$e$ are optimal it follows that $c(\binom{B_i}{2})=0$, for $i\in \Sigma_0$ 
and 
$c(B_i\times B_j)=0$, for $(i,j)\in \Sigma_c$ (here, for a subset $E$ of 
edges, $c(E) = \sum_{(i,j)\in E} (c_i + c_j)$). Since $d$ and $e$ are the 
only optimal solutions it now follows that
$$\delta(\Sigma_0) + \#\Sigma_c = 1,$$
where $\delta(\Sigma_0) = 1$ if $\#\Sigma_0 = 1$ and
the unique element $i$  of $\Sigma_0$ 
satisfies $\#B_i \geq 2$ and $\delta(\Sigma_0)=0$ otherwise. The result 
follows. $\Box$

In order to determine the higher dimensional faces of 
$\TPn$ we need to 
know more about the interaction between the lattice structure of $TP(n)$ 
and the faces of $\TPn$. In particular, we would like to know 
the 
answer to the following question: given a face of $\TPn$ the set of 
elements of $TP(n)$ lying on this face is closed under $\wedge$ and 
$\vee$. What are the join irreducible elements on this face?.

\bt
For $n\geq 3$, the number of edges of $\DPn$
is $2^{n-2}(2n-3)$.
\et
\pf Let $TD(n)$ denote the set of all threshold partitions $d=(d_1,\ldots 
,d_n)\in TP(n)$ satisfying $d_1 = n-1$ (or, equivalently, $d_n \not= 0$). 
Given $d\in TD(n)$ as above, let $m(d)$ denote the number of dominating 
vertices in $d$, i.e., $m(d)$ is the largest index $j$ with $d_j = n-1$. 
We assert that
\beq \label{ndp}
\sum_{d\in TD(n)} m(d) = 2^{n-1}. \eeq
The proof is by induction on $n$, the cases $n=1,2$ being checked easily. 
Let $n\geq 3$ and let $d=(d_1,\ldots ,d_n)\in TD(n)$. Write
$$
(d_1,\ldots ,d_n) = (n-1,1,\ldots ,1) + (0,e_1,\ldots ,e_{n-1}),$$
where $e=(e_1,\ldots ,e_{n-1})\in TP(n-1)$.

If $e_{n-1}=0$, then $m(d)=1$. The number of $e\in TP(n-1)$ with 
$e_{n-1}=0$ 
is $2^{n-3}$. If $e_{n-1}\not=0$, then $m(d)=m(e)+1$.
The number of $e\in TP(n-1)$ with $e_{n-1}\not=0$ 
is also $2^{n-3}$. 

It follows by induction that
$$\sum_{d\in TD(n)} m(d) = 2^{n-3} + 2^{n-3} + \sum_{e\in TD(n-1)} m(e) = 
2^{n-3} + 2^{n-3} + 2^{n-2} = 2^{n-1}.$$

Let $E_n$ denote the number of edges of $\DPn$. We prove the following 
recurrence
$$E_n = 2E_{n-1} + 2^{n-1},\;\;n\geq4,$$
with $E_3 = 6$. The result follows easily from this recurrence by 
induction. 

The polytope $\mathsf{DP(3)}$ has 4 vertices any two of which are 
adjacent, so $E_3 = 6$. Let $n\geq 4$ and let $d=(d_1,\ldots ,d_n)$ 
and $e=(e_1,\ldots ,e_n)$ be adjacent extreme points of $\DPn$.  We say 
that $d$ and $e$ are {\em straight neighbors} if both $d$ and $e$ are 
dominating (i.e., $d_1 = e_1 = n-1$) or both $d$ and $e$ are isolated 
(i.e., $d_n = e_n = 0$) and we say that $d$ and $e$ are {\em cross 
neighbors} if one of them is dominating and the other isolated. Straight 
neighbors lie on one of the facets $x_1 \leq n-1$ or $x_n \geq 0$ and 
thus, by induction, the number of straight neighbors is $2E_{n-1}$. We 
shall now show that the number of cross neighbors is $2^{n-1}$. This will 
prove the recurrence.

Let $e=(e_1,\ldots ,e_n)\in TP(n)$ with $e_n = 0$. We want to count the 
number of cross neighbors of $e$. The following two cases arise:

(i) $e=(0,\ldots ,0)$: If $d$ is a cross neighbor of $e$ then $d-e$ will 
have nonzero first and last components and thus, by Theorem \ref{edges}, 
$d= \Omega_{[n]}$ or $\Omega_{I,J}$, with $1\in I$ and $n\in J$. It is 
easily checked that the only possibility for $(I,J)$ is $I=\{1\}$ and 
$J=\{2,\ldots ,n\}$. Thus $e$ has two cross neighbors.

(ii) $e\not= (0,\ldots ,0)$: Let $e_1 = j$, where $1\leq j \leq n-2$. Then 
$e_{j+2}=\cdots =e_n=0$. Put $e'=(e_1,\ldots,e_{j+1})\in TD(j+1)$. 
A small calculation  using Theorem \ref{edges} shows that $d\in TP(n)$ is 
a cross neighbor of $e$ if and only if
$$ d-e = \Omega_{I,\{j+2,\ldots ,n\}},$$
where $I=\{1,\ldots ,k\}$ for some $k\leq m(e')$. Thus $e$ has $m(e')$ 
cross neighbors.

Using (i), (ii) above and the formula (\ref{ndp}) we see that the total 
number of cross neighbors is equal to
$$2 + \sum_{j=2}^{n-1} 2^{j-1} = 2 + 2(2^{n-2} - 1) = 2^{n-1}.$$

That completes the proof. $\Box$

\begin{center}\section{{\large  Majorization and $r$-ideal partitions}}  
\end{center}

Let $a=(a(1),\ldots ,a(n))$ and $b=(b(1),\ldots ,b(n))$ be real sequences 
of 
length $n$. Denote the $i$-th largest component of $a$ (respectively, $b$) 
by $a[i]$ (respectively, $b[i]$). We say that $a$ {\em majorizes} 
$b$, denoted $a \succeq b$, if
$$ \sum_{i=1}^k a[i] \geq \sum_{i=1}^k b[i],\;\;k=1,\ldots ,n,$$
with equality for $k=n$. The majorization is {\em strict}, denoted $a 
\succ b$, if at least one of the inequalities is strict, namely if $a$ is 
not a permutation of $b$.

We recall a fundamental lemma about majorization in integer 
sequences, 
called Muirhead's lemma. If $a=(a(1),\ldots ,a(n))$ is a sequence and 
$a(i)\geq a(j) + 2$, then the following operation is called a {\em unit 
transformation} from $i$ to $j$ on $a$: subtract 1 from $a(i)$ and add 1 
to $a(j)$. Clearly, if $b$ is obtained from $a$ by a sequence of unit 
transformations, then $a\succ b$. The converse is also true for integer 
sequences.

\bt
\label{mhl}
(Muirhead Lemma) If $a$ and $b$ are integer sequences and $a\succeq b$, 
then some permutation of $b$ can be obtained from $a$ by a finite sequence 
of unit transformations.
\et
For a proof see {\bf \cite{mp}, \cite{mo}}.

\bl
\label{gsl}
Let $a$ be a $r$-graphical sequence (i.e., the degree sequence of an
$r$-graph) and let $b$ be a nonnegative integral 
sequence such that $a\succeq b$. Suppose that $a=a_0,a_1,\ldots ,a_l=b$ 
are nonnegative integral sequences such that, for $1\leq i\leq l-1$, $a_i$ 
yields $a_{i+1}$ by a 
unit transformation. Then each 
$a_i$, in particular $b$, is an $r$-graphical sequence.
\el
\pf By induction on $i$, $a_0$ being $r$-graphical by hypothesis. Let 
$H_i$ be an $r$-graph with degree sequence $a_i$ and let $a_{i+1}$ 
be obtained by a unit transformation from $j$ to $k$ on $a_i$, 
where $a_i(j)\geq a_i(k) + 2$. Since 
the vertex $j$ has larger degree than the vertex $k$ in $H_i$, there must be a 
$X\in S(n,r-1)$ with $j,k \not\in X$
such that $X\cup \{j\}$ is an edge of $H_i$ and 
$X\cup \{k\}$ is not an edge of $H_i$. The simple $r$-graph $H_{i+1}$ 
obtained from $H_i$ by deleting the edge $X\cup \{j\}$ and adding the edge 
$X\cup \{k\}$ has the degree sequence $a_{i+1}$.$\Box$ 

\vskip 1ex
We can now prove our result characterizing $r$-graphical partitions.

\noi
\pf (of Theorem \ref{mhd}) (only if) Let $d=(d(1),\ldots ,d(n))$ be the 
degree sequence of the 
$r$-graph $H = 
([n], E)$. 

Suppose that  there are distinct vertices $i,j\in [n]$,  
$X\in 
S(n,r-1), i,j\not\in X$ such that $d(i)\geq d(j)$, $X\cup \{i\} 
\not\in E$, and 
$X\cup\{j\} \in E$. Let $H'$ be the $r$-graph obtained from $H$ by 
deleting the edge $X\cup\{j\}$ and adding $X\cup\{i\}$ as an edge. Then 
the  
degree sequence $d'$ of $H'$ strictly majorizes $d$ and $d$ can be 
obtained from $d'$ by a unit transformation. We say that $H'$ is 
obtained from $H$ by a reverse unit transformation.

Since the number of nonnegative integer sequences of length $n$  with 
sum $\sum_{i=1}^n d(i)$
is finite and every application of a reverse unit transformation makes the 
degree sequence go strictly up (in the majorization order)
we may, starting from $H$, apply reverse unit 
transformations a finite 
number of times (possibly zero) to obtain a hypergraph $H'' = ([n], E'')$ 
whose degree sequence $d''$ majorizes $d$ and such that we cannot apply 
any reverse unit transformation to $H''$. That is, for all $i,j\in 
[n],i\not=j$ and $X\in E$ 
$$d''(i)\geq d''(j) \mbox{ and } j\in X \mbox{ implies } i\in X \mbox{ or 
} (X-\{j\})\cup \{i\} \in E.$$
It follows that $H''$ is an order ideal of $S(n,r)$ under the order 
induced from the  linear order 
on $[n]$ 
obtained by listing the vertices in weakly decreasing order of their 
degrees (in $H''$). The result follows.

\noi (if) This follows from Muirhead's Lemma and Lemma \ref{gsl}. $\Box$

\end{document}